\documentclass[12pt]{article}
\usepackage{authblk}
\usepackage[english]{babel}
\usepackage[utf8x]{inputenc}
\usepackage[colorinlistoftodos]{todonotes}
\usepackage{graphicx, fleqn, tabularx, wrapfig}
\usepackage{amsthm,amsmath,textcomp}
\usepackage[T1]{fontenc}
\usepackage{charter}
\usepackage[bitstream-charter]{mathdesign}
\usepackage{indentfirst}
\usepackage{caption}
\usepackage{subcaption}
\usepackage{booktabs}
\usepackage{hyperref}
\usepackage{longtable}

\newtheorem{theorem}{Theorem}[section]
\newtheorem{lemma}[theorem]{Lemma}

\title{\Large \textbf{A Note on Weighted Fermat Problem in $\mathcal{L}^p$ Space}}
\author{Shikun Liu}
\date{}

\begin{document}
\maketitle

\begin{abstract}
In this paper, we begin by introducing a well-known geometry concept: the Fermat point in a triangle. Then, we generalize the problem and propose an iterative algorithm based on gradient descent to the weighted form in $\mathcal{L}^p$ space. We also build specific solutions to some special norms: one, two and infinity. We show that the solution may not be unique in norm-one and infinity. Finally, we provide our qualitative results in the PASCAL-based program.

\end{abstract}

\section{Introduction}

The Fermat point optimization problem was named after the French mathematician \textit{Pierre de Fermat} in a private letter in 17th century: Given three points in the plane, finding the fourth point such that the sum of the distances to the three given points is the minimum. The problem was generalized by \textit{Simpson}$^{[1]}$ to ask finding the minimum weighted sum of distances from three given points.

In this note, we consider the problem in $\mathcal{L}^p$ space. Let ${\bf y}_i=(y^i_1,y^i_2,\cdots,y^i_n)$ be fixed points with their positive weights $k_i$ for $i=1,2,\cdots, m$ for total $m$ points. Given the optimized point ${\bf F}=(x_1,x_2,\cdots,x_n)$ in $\mathbb{R}^n$, we want to find $(x_1,x_2,\cdots,x_n)$ such that
$$
\sum_{i=1}^m\left(\sum_{j=1}^n k_i|y_j^i-x_j|^p  \right)^{1/p}\!=\text{min}.
$$

\section{Fermat Point in Manhattan Norm}
In Manhattan norm (norm-1 space), the distance equation is defined as:
$$
\mathcal{L}^1=||{\bf x}-{\bf y}||_1=\sum_{i=1}^{n}|x_i-y_i|.
$$

According to the definition of Manhattan norm, we propose our target function,
$$
f(x_1,x_2,\cdots, x_n)=\sum_{i=1}^m\sum_{j=1}^n k_i|y_j^i-x_j|.
$$

Reorder coordinate sets $\left\{y_j^i\right\}$ in each $j^{th}$ dimension an increasing order $\left\{\widetilde{y}_j^i\right\}$ with its weight $\widetilde{k}^{i}_j$. Define,
$$
\widetilde{f}_j(x_j)=\sum_{i=1}^m \widetilde{k}^i_j |\widetilde{y}_j^i-x_j|,\quad j=1,2,\cdots,n
$$

Since the function $\widetilde{f}_j(x_j)$ are independent in each dimension, the optimal point can be found when the derivative in each dimension of $\widetilde{f}_j(x_j)$ is zero or NOT exist.

Let $\widetilde{m}_{l_j}^{i}<x_j<\widetilde{m}_{l_j+1}^{i}, l_j=1,2,\cdots,m-1$. Then, $\widetilde{f}_j(x_j)$ can be split into,
$$\widetilde{f}_j(x_j)=\sum_{i=1}^{l_j} \widetilde{k}^i_j |\widetilde{y}_j^i-x_j|+\sum_{i=l_j+1}^{m} \widetilde{k}^i_j|x_j-\widetilde{y}_j^i|.$$

Therefore, the optimal point ${\bf F}(x_1,x_2,\cdots,x_n)$ holds when the function satisfies,
$$
 \frac{\partial\widetilde{f}_j^-}{\partial x_j^-}= \sum_{i=1}^{l_j}\widetilde{k}_j^i-\sum_{i=l_j+1}^{m}\widetilde{k}_j^i \leq 0, \quad \frac{\partial f_1^+}{\partial x_j} =  \sum_{i=1}^{l_j}\widetilde{k}_j^i-\sum_{i=l_j+1}^{m}\widetilde{k}_j^i \geq 0,\quad j = 1,2,\cdots,n
$$

From the solutions above, we can see that Fermat point in norm-1 space may have infinity solutions. The visual representations with random 50 points in $\mathbb{R}^2$ space can be illustrated as follow.

\begin{figure}[ht!]
\centering
\includegraphics[scale=0.33]{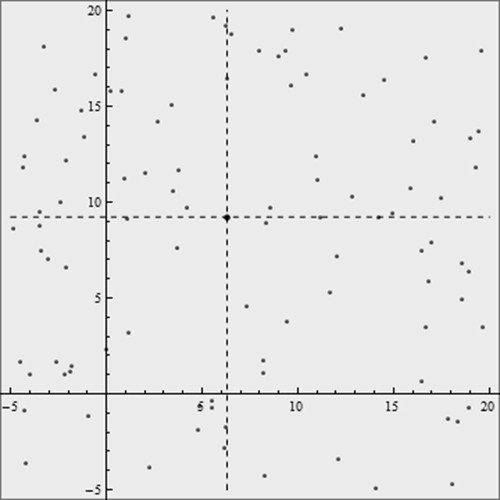}
\qquad \qquad \qquad
\includegraphics[scale=0.33]{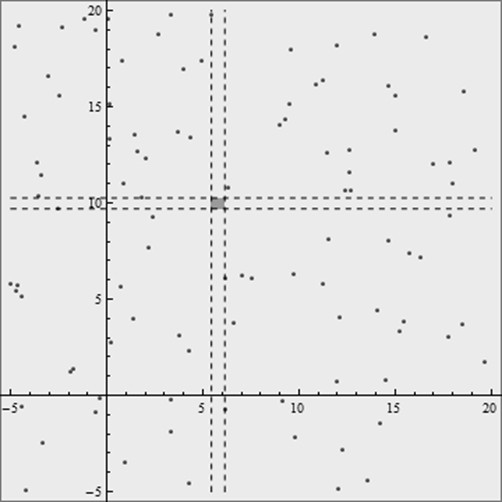}
\caption{Two possible cases of Fermat point(s) (marked in grey point or box) in $\mathcal{L}^1$ norm, in which left has one unique solution and right has infinity number of solutions. }
\end{figure}

\section{Fermat Point in Euclidean Norm}
The unweighted 3-points problem in the Euclidean Norm is the original Fermat Problem. A solution to the one is to construct equilateral triangles on the sides with the vertices pointing outward. The intersected vertex from three circles is the Fermat point. The visualization is Figure 2, in which the black points are fixed points, and the grey point is the Fermat Point.
\begin{figure}[htbp]
\includegraphics[scale=0.23]{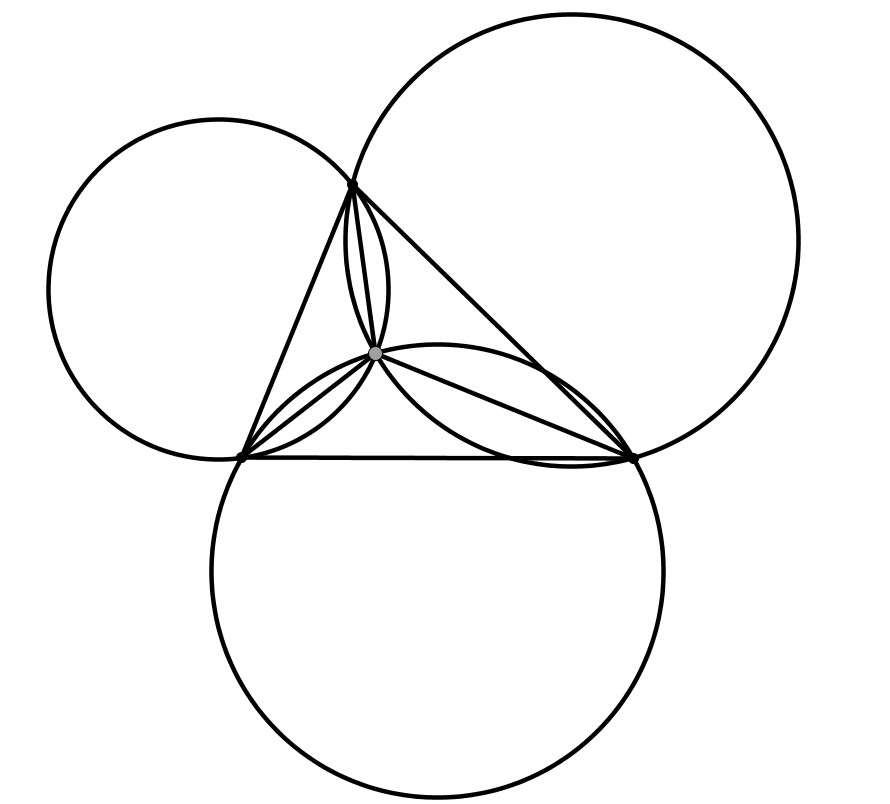}
\centering
\caption{Fermat point in a triangle}
\end{figure}

\subsection{The Varignon Frame}
A good tool that allow us to inspire the general Fermat problem in Euclidean space is a mechanical analogue device \textit{The Varignon Frame} named after the French mathematician \textit{Varignon}. A board is drilled with $m$ holes corresponding to the coordinates of the $m$ fixed points. One string is fed through each hole and all strings are tied together in one knot. Below the board, weights are attached to the strings and proportional to the board. Given only gravity we concern and we ignore friction, the knot will settle at the optimal point. A graph of the Varignon frame is shown in Figure 3.$^{[2]}$

\begin{figure}[h]
\centering
\includegraphics[scale=0.3]{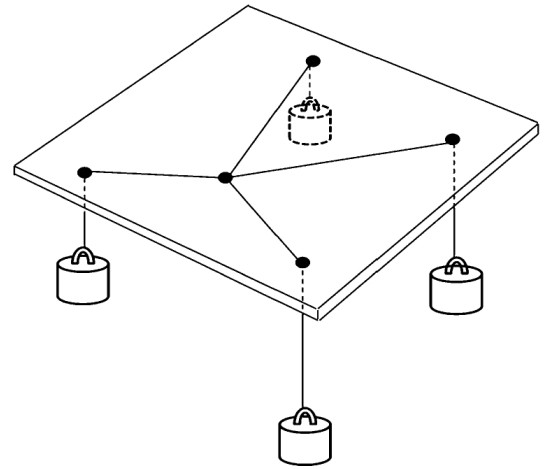}
\caption{The Varignon Frame}
\end{figure}

\begin{lemma}
The minimal potential energy point is the optimal point.
\end{lemma}
\begin{proof}
Let the zero potential energy level is the ground, and height of this system is $h$, the length of each string is $l_i$, and its weight is $G_i$, where $i=1,2,\cdots, m$. We set up an initial point $P$ in this system, then we have
$$
E(P)=\sum_{i=1}^m G_i [h-\left(l_i-\overline{PM_i}\, \right)]=\sum_{i=1}^m G_i(h-l_i)+\sum_{i=1}^m G_i (\, \overline{P M_i}\, )
$$

When $P$ is the minimal potential energy point, that is $E(P)$ is the minimum,  $\sum_{i=1}^m G_i (\, \overline{P M_i}\, )$ reaches the minimum as $ \sum_{i=1}^m G_i(h-l_i)$ is the constant.
\end{proof}

\subsection{The Weiszfeld's Algorithm}
The Weiszfeld's iterative algorithm$^{[3]}$ was proposed to the Fermat problem in $\mathbb{R}^2$ space. The idea behind this algorithm, starting anywhere and trying to converge to ${\bf F}$, is actually very simple. Let
$$
{\bf F}=\sum_{i=1}^m\left(\sum_{j=1}^n k_i|y_j^i-F|^2 \right)^{1/2}.
$$

We assume ${\bf F} \notin {\bf y}_i$, then the negative of the gradient of $f$ at ${\bf F}$ equals,
$$
R({\bf F})=\sum_{i=1}^m k_i \frac{{\bf y}_i-{\bf F}}{||{\bf y}_i-{\bf F}||}.
$$

Therefore, let  $R({\bf F})=0$ which is equivalent to,
$$T({\bf F})= \sum_{i=1}^m \frac{k_i {\bf y}_i}{\|{\bf F}-{\bf y}_i\|}/{\sum_{i=1}^m \frac{k_i }{\|{\bf F}-{\bf y}_i\|}},$$
as our iterative algorithm.

\textit{Rautenbach}$^{[7]}$ proposed an extension $T^*$ of $T$ to $\mathbb{R}^n$ space. Let $P \in \mathbb{R}^n,$
$$
\begin{cases}
\text{If }T^*(P)=P,  \text{ if and only if }P={\bf F}& (i) \\
\text{If } T^*(P)\neq P,\text{ then } f(T^*(P))<f(P)& (ii)
\end{cases}
.
$$

\textit{Kuhn}$^{[4]}$ proved the part $(ii)$  and showed that, given an initial point $X_0 \in \mathbb{R}^n$, the sequence of points $\left\{X_n \right\}_{n\in N}$ generated by the iterative process $X_{k+1}=T\left(X_k\right)$ converges to the optimal solution of the Fermat problem if no point in the sequence is the vertex. In this last case, Kuhn claimed to have proved that this fact only occurs for a denumerable set of initial points. This result was modified by \textit{Chandrasekaran} and \textit{Tamir} $^{[5]}$, who gave two counterexamples and established the following conjecture: \textit{"If the convex hull of the set of vertices is of full dimension, then the set of initial
points for which the sequence generated by the Weiszfeld's algorithm yields in a vertex is
denumerable."}

\textit{Vardi} and \textit{Zhang}$^{[6]}$ improved the limit of Weiszfeld's algorithm in case of ${\bf F} \in {\bf y}_i$ appears in any iteration.

\section{Fermat Point in Chebyshev Norm}
\quad \, The two dimensional Chebyshev distance has the similar property as the Manhattan distance. Hence, it can easily be viewed as the transformation of Manhattan distance, with the transformation matrix is
$$
T = \begin{bmatrix}
       \frac{\sqrt{2}}{2} & 0 \\[0.3em]
       0& \frac{\sqrt{2}}{2}       \\[0.3em]
      \end{bmatrix}
 \begin{bmatrix}
      \cos{\frac{\pi}{2}} &- \sin{\frac{\pi}{2} }  \\[0.3em]
       \sin{\frac{\pi}{2}}& \quad \! \cos{\frac{\pi}{2}}       \\[0.3em]
      \end{bmatrix}
.
$$
\begin{lemma}
$\mathcal{L}^\infty ({\bf y}^i,{\bf F})=\mathcal{L}^1 (T({\bf y}^i),T({\bf F}))$
\end{lemma}

\begin{proof}

\begin{equation}\label{e:barwq}\begin{split}
\mathcal{L}^1 (T(y_1^i),T(F)) & = \frac{\sqrt{2}}{2}\left| y_1^i \cos \frac{\pi}{2}+y_2^i \sin \frac{\pi}{2}-x_1\sin \frac{\pi}{2}-x_2 \sin \frac{\pi}{2}\right|  \\
&+\frac{\sqrt{2}}{2}\left| -y_1^i \cos \frac{\pi}{2}+y_2^i \sin \frac{\pi}{2}+x_1\sin \frac{\pi}{2}-x_2 \sin \frac{\pi}{2}\right|   \\
&=\frac{\sqrt{2}}{2} \max \bigg\{\frac{1}{\sqrt{2}} \max \big\{(y_2^i-x_2+y_2^i-x_2),(y_1^i-x_1+y_1^i-x_1), \nonumber \\
&(-y_1^i+x_1-y_1^i+x_1),(-y_2^i+x_2-y_2^i+x_2) \big\} \bigg\} \nonumber \\
&=\max\left\{|y_1^i-x_1|,|y^i_2-x_2|\right\}.
\end{split}\end{equation}

\end{proof}

However, the equivalence between $\mathcal{L}^1 \text{ and } \mathcal{L^{\infty}}$ cannot be generalized to the $\mathbb{R}^n$ due to the different properties on these two norms in higher dimensions.

\section{Fermat Point in $P$-Norm}
The $\mathcal{L}^p$ Space are the certain vector spaces of measurable functions. A vector space on which a norm is defined is called a normed vector space. It satisfies the following conditions,$^{[8]}$
\begin{itemize}
\item $||f||>0$ for $f\neq 0$
\item $||0||=0$
\item $||a\cdot f||=|a|\cdot ||f||$ for all $a\in \mathbb{C}^1$
\item$ ||f+g||\leq ||f||+||g||$ The Triangle Inequality.
\end{itemize}

For $0<p\leq \infty$ we will denote the set of all $\mu$-measurable functions, $f$, such that $||f||^p$ is integrable by $\mathcal{L}^p (X,\mathcal{F},\mu)$. We define a real-valued function, $||\cdot||_p$, on $\mathcal{L}^p$ by ,
$$
||f||_p= \left(\int ||f(x)||^p\,  d\mu (x)  \right)^{1/p}.
$$

For $1\leq p<\infty$ this function is called the $\mathcal{L}^p$-norm. For $0<p<1$, $\mathcal{L}^p$ is still a vector space, but $||\cdot||_p$ is no longer a norm. Although $||\cdot||_p^p$ can be used to define an interesting metric on $\mathcal{L}^p$, we will not consider this case further.
\begin{lemma}
If the points are not collinear, then $f$ is strictly convex.
\end{lemma}
\begin{proof}
For $X,Y,Q\in \mathbb{R}^n$, and $0<t<1$, the Minkowski inequality shows,
$$ ||X+Y||_p\leq ||X||_p+||Y||_p$$

implies,
\begin{align}
||t(X-&Q)+(1-t)(Y-Q)||^p \nonumber \\
 &\leq  ||t(X-Q)||^p + 2||t(X-Q)||\cdot ||(1-t)(Y-Q)||+||(1-t)(Y-Q)||^p \nonumber \\
&=(t||X-Q||+(1-t)||Y-Q||)^p \nonumber
\end{align}

Hold if and only if $X,Y,Q$ are not collinear, therefore,
$$f(tX+(1-t)Y)<tf(X)+(1-t)f(Y)$$
\end{proof}

Since $f(x_1,x_2,\cdots, x_n)$ is convex, it proves that the Fermat point exists.  Thus, by gradient descent
$$\nabla f=\frac{\partial f}{\partial x_l}(x_1,x_2,\cdots,x_n=0,l=1,2,\cdots,n,
$$ the function gets the minimum. It gives,
$$
\frac{\partial f}{\partial x_l}(x_1,x_2,\cdots,x_n)=\sum_{i=1}^m \sum _{j=1}^n k_i |y_l^i-x_l|^{p-1}\left( \left|  m_j^i-x_j)\right|^p  \right)^{1/p-1} \frac{\sqrt{y_l^i-x_l}}{y_l^i-x_1}=0,l=1,2,\cdots,n.
$$

Therefore, we have our generalized iterative algorithm,
$$
T(x_l)=\frac{\sum_{i=1}^m \sum _{j=1}^n k_i x_l^i |y_l^i-x_l|^{p-1}\left( \left|  y_l^i-x_j)\right|^2  \right)^{p/2-1}}{\sum_{i=1}^m \sum _{j=1}^n k_i |y_l^i-x_l|^{p-1}\left( \left|  y_l^i-x_j)\right|^2  \right)^{p/2-1}}, \quad l=1,2,\cdots,n.
$$

\section{The Implement Test}
\quad \, We build the PASCAL-based programs for iterative algorithms we proposed in the above sections. The following tables below are the test results for random test points. In norm-1 or $\infty$, the solution gives the exact result or range; in norm-2 and $p$, the programs also add the precision as a condition. We set the first iteration as the gravity point of the system.

\smallskip

\begin{table}[ht!]
\centering
\footnotesize
\begin{tabular}{@{}lccc@{}}
\toprule
\multicolumn{1}{c}{\textbf{Input:}} &\quad  \quad Number of Points & \quad \quad \quad Dimension & Weights \\
\multicolumn{1}{c}{ } & \quad \quad 4 & \quad \quad \quad 6 & 2, 3, 1, 2\\ \midrule
\multicolumn{4}{l}{\textbf{Test Points:} (4, 11, 3, 4, 5, 6), (13, 13, 2, 1, 5, 7), (17, 6, 8, 6, 7, 8), (8, 6, 4, 6, 7, 9)}  \\
\multicolumn{4}{l}{\textbf{Output:} \quad \quad(8$\sim$13, 11, 3, 4, 5, 7)   }          \\
\bottomrule
\end{tabular}
\caption{$\mathcal{L}^{p} $ Algorithm}
\end{table}

\newpage

\begin{table}[ht!]
\centering
\footnotesize
\begin{tabular}{@{}lcccc@{}}
\toprule
\multicolumn{1}{l}{\textbf{Input:}} & \quad Number of Points  & \quad  \quad  Precision &  \quad  \quad  \quad Dimension & Weights \\
\multicolumn{1}{c}{ } & \quad 6  & \quad  \quad 0.001& \quad  \quad  \quad 2 &  1, 2, 4, 7, 6, 5\\ \midrule
\multicolumn{5}{l}{\textbf{Test Points:} (4.71, -1.84), (-3.15, -2.44), (0.17, 2.99),(6.35, 2.86), (5.55, 2.44), (3.22, -2.56)}  \\
\multicolumn{5}{l}{\textbf{Output:} Start : (1.480, -0.140) } \\
\multicolumn{5}{l}{\quad \quad \quad \quad \quad \quad \! 1 : (1.713, 0.086)} \\
\multicolumn{5}{l}{\quad \quad \quad \quad \quad \quad \! 2	: (1.763, 0.124)} \\
\multicolumn{5}{l}{\quad \quad \quad \quad \quad \quad \! 3	: (1.776, 0.126)}\\
\multicolumn{5}{l}{\quad \quad \quad \quad \quad \quad \! 4	: (1.782, 0.122)}\\
\multicolumn{5}{l}{\quad \quad \quad \quad \quad \quad \! 5	: (1.787, 0.118)}\\
\multicolumn{5}{l}{\quad \quad \quad \quad \quad \quad \! 6	: (1.790, 0.115)}\\
\multicolumn{5}{l}{\quad \quad \quad \quad \quad \quad \! 7	: (1.793, 0.112)}\\
\multicolumn{5}{l}{\quad \quad \quad \quad \quad \quad \! 8 :  (1.795, 0.110)}\\
\multicolumn{5}{l}{\quad \quad \quad \quad \quad \quad \! 9	: (1.798, 0.108)}\\
\multicolumn{5}{l}{\quad \quad \quad \quad \quad \quad 10: (1.799, 0.106)}\\
\multicolumn{5}{l}{\quad \quad \quad \quad \quad \quad 11: (1.801, 0.105)}\\
\multicolumn{5}{l}{\quad \quad \quad \quad \quad \quad 12: (1.802, 0.104)}\\
\bottomrule
\end{tabular}
\caption{$\mathcal{L}^{2} $ Algorithm}
\end{table}

\begin{table}[ht!]
\centering
\footnotesize
\begin{tabular}{@{}lccccc@{}}
\toprule
\multicolumn{1}{l}{\textbf{Input:}} & Number of Points  & Norm &\quad    Precision &  \quad Dimension & Weights \\
\multicolumn{1}{c}{ } & 5  & 2.4 & \quad 0.00001 & \quad 5 &5, 9, 1, 8, 6\\ \midrule
\multicolumn{6}{l}{\textbf{Test Points:} (8, 5, 4, 8, 3), (3, 3, 7, 6, 3), (8, 7, 2, 6, 6), (4, 9, 3, 6, 2), (5, 6, 4, 5, 4)}  \\
\multicolumn{6}{l}{\textbf{Output:} Start : (4.72414, 5.75862, 4.58621, 6.13793, 3.03448) } \\
\multicolumn{6}{l}{\quad \quad \quad \quad \quad \quad \! 1	: (4.95849, 5.91367, 4.57670, 5.74708, 3.46078)} \\
\multicolumn{6}{l}{\quad \quad \quad \quad \quad \quad \! 2	: (4.93330, 5.93940, 4.40332, 5.61458, 3.46066)} \\
\multicolumn{6}{l}{\quad \quad \quad \quad \quad \quad \! 3	: (4.94537, 5.99067, 4.38619, 5.56072, 3.53263)}\\
\multicolumn{6}{l}{\quad \quad \quad \quad \quad \quad \! 4	: (4.94590, 5.99950, 4.35138, 5.52565, 3.54720)}\\
\multicolumn{6}{l}{\quad \quad \quad \quad \quad \quad \! 5	: (4.94774, 6.01196, 4.34204, 5.50728, 3.56604)}\\
\multicolumn{6}{l}{\quad \quad \quad \quad \quad \quad \! 6	: (4.94833, 6.01319, 4.33222, 5.49539, 3.57347)}\\
\multicolumn{6}{l}{\quad \quad \quad \quad \quad \quad \! 7	: (4.94886, 6.01513, 4.32840, 5.48840, 3.57989)}\\
\multicolumn{6}{l}{\quad \quad \quad \quad \quad \quad \! 8	: (4.94909, 6.01562, 4.32523, 5.48396, 3.58315)}\\
\multicolumn{6}{l}{\quad \quad \quad \quad \quad \quad \! 9	: (4.94928, 6.01613, 4.32359, 5.48123, 3.58549)}\\
\multicolumn{6}{l}{\quad \quad \quad \quad \quad \quad 10: (4.94938, 6.01634, 4.32243, 5.47950, 3.58684)}\\
\multicolumn{6}{l}{\quad \quad \quad \quad \quad \quad 11: (4.94945, 6.01651, 4.32176, 5.47842, 3.58773)}\\
\multicolumn{6}{l}{\quad \quad \quad \quad \quad \quad 12: (4.94949, 6.01660, 4.32131, 5.47774, 3.58828)} \\
\multicolumn{6}{l}{\quad \quad \quad \quad \quad \quad 13: (4.94952, 6.01666, 4.32104, 5.47730, 3.58864)} \\
\multicolumn{6}{l}{\quad \quad \quad \quad \quad \quad 14: (4.94953,  6.01669, 4.32086, 5.47703, 3.58886)} \\
\multicolumn{6}{l}{\quad \quad \quad \quad \quad \quad 15: (4.94954, 6.01672, 4.32075, 5.47686, 3.58900)} \\
\multicolumn{6}{l}{\quad \quad \quad \quad \quad \quad 16: (4.94955, 6.01673, 4.32068, 5.47675, 3.58909)} \\
\multicolumn{6}{l}{\quad \quad \quad \quad \quad \quad 17: (4.94956, 6.01674, 4.32063, 5.47668, 3.58914)} \\
\multicolumn{6}{l}{\quad \quad \quad \quad \quad \quad 18: (4.94956, 6.01675, 4.32061, 5.47664, 3.58918)} \\
\multicolumn{6}{l}{\quad \quad \quad \quad \quad \quad 19: (4.94956, 6.01675, 4.32059, 5.47660, 3.58920)} \\
\multicolumn{6}{l}{\quad \quad \quad \quad \quad \quad 20: (4.94956, 6.01675, 4.32058, 5.47659, 3.58921)} \\
\bottomrule
\end{tabular}
\caption{$\mathcal{L}^{p} $ Algorithm}
\end{table}

\newpage
\section*{Acknowledgements}
I would like to express my great gratitude to my adviser {\em Donglei Gu }and my friend {\em Wilhelm Braum}, for giving many insightful advices in this paper. Furthermore, I would like to thank my friend {\em Liaoyuan Shen} for the technical help in programming.

\end{document}